\documentclass[11 pt]{amsart}
\usepackage{amsmath,amsthm,amsfonts,amssymb,amscd,amsrefs}
\usepackage[top=3cm, bottom=3cm,left=2cm, right=2cm]{geometry}
\allowdisplaybreaks[1]
\usepackage{enumitem}
\setlist[enumerate]{font=\normalfont}
\usepackage{centernot}
\usepackage{xcolor}
\usepackage{soul}
\usepackage{tikz}
  \usepackage[linkcolor=red,citecolor=green,colorlinks=true]{hyperref}
\usepackage{mathrsfs}
\usepackage{mathtools}
\usepackage{etoolbox}
\usepackage{scalerel,stackengine}
\stackMath
\newcommand\widecheck[1]{%
\savestack{\tmpbox}{\stretchto{%
  \scaleto{%
    \scalerel*[\widthof{\ensuremath{#1}}]{\kern-.6pt\bigwedge\kern-.6pt}%
    {\rule[-\textheight/2]{1ex}{\textheight}}
  }{\textheight}%
}{0.5ex}}%
\stackon[1pt]{#1}{\scalebox{-1}{\tmpbox}}%
}
\theoremstyle{definition}
\newtheorem{defn}{Definition}[section]

\theoremstyle{remark}

\newtheorem{rmk}[defn]{Remark}
\theoremstyle{plain}
\newtheorem{thm}[defn]{Theorem}
\newtheorem{lem}[defn]{Lemma}
\newtheorem{cor}[defn]{Corollary}

\DeclareMathOperator*{\ess}{ess\,sup}

\numberwithin{equation}{section}
\makeatletter 

\newcommand{\Rmnum}[1]{\expandafter\@slowromancap\romannumeral #1@}
\newcommand{\sgn}{\mathop{\mathrm{sgn}}}
\makeatother
\makeatletter
\def\imod#1{\allowbreak\mkern10mu({\operator@font mod}\,\,#1)}
\makeatother
\makeatletter 
\patchcmd{\@settitle}{\uppercasenonmath\@title}{}{}{}
\patchcmd{\@setauthors}{\MakeUppercase}{}{}{}
\patchcmd{\section}{\scshape}{}{}{}
\makeatother
\keywords{Strong subdifferentiability, Bochner $L^{p}$-space, Radon measure, Separability}
\subjclass[2020]{ Primary 46E40, 46B20, secondary 28B20}
\begin{document}
\selectfont
\title{Strong Subdifferentiability in Bochner Function spaces  $L^{p}(\mu, X)$}
 \author{Saurabh Dwivedi}
 \author{C. R. Jayanarayanan}
     \address{Department of Mathematics, Shiv Nadar Institution of Eminence, Delhi-NCR, India}
	 \email{sd605@snu.edu.in}
	 \address{Department of Mathematics, Indian Institute of Technology Palakkad, 678557, India}
	 \curraddr{Department of Mathematics, Indian Institute of Technology Palakkad, 678557, India}
	 \email{crjayan@iitpkd.ac.in}
	
	\begin{abstract}
    In this article, for $1\leq p<\infty$, we provide a characterization of strongly subdifferentiable points of  $L^{p}(\mu, X)$, when $X^*$ is separable and $\mu$ is the Lebesgue measure on $[0, 1]$. For  $f\in L^{p}(\mu, X)$ with  $\|f\|_{p}=1$ we show that the norm of $L^{p}(\mu, X)$ is strongly subdifferetiable at $f$ if and only if the norm of $X$ is strongly subdifferentiable at $f(t)$ for almost every \(t\in[0,1]\) with $f(t)\neq 0$. Moreover, for a separable Banach space $X$, we show that the norm of $L^{p}(\mu, X)$ is strongly subdifferentiable if and only if the norm of $X$ is strongly subdifferentiable.
	\end{abstract}
    \maketitle
\section{Introduction}

Lebesgue-Bochner spaces $L^p(\mu, X)$
constitute a natural vector-valued extension of the classical $L^p$-spaces, and the geometric properties of these spaces, particularly the G\^ateaux and Fr\'echet  differentiability of their norms, have been extensively studied in the literature (see \cites{LSK, LSK1, LSK2, WN}). For $1<p<\infty$, G\^ateaux differentiability of the norm of $L^p(\mu, X)$ was characterized in \cite{LSK}, while Wang and Nan, in \cite{WN},  characterized its points of Fr\'echet differentiability in terms of
the Fr\'echet differentiability of their pointwise values
 in the underlying Banach space $X$.  Strong subdifferentiability  of the norm, introduced in \cite{DG}, is a nonsmooth extension of  Fr\'echet differentiability of the norm. It is closely related to the continuity properties of the  duality mapping and has significant connections with several geometric and approximation-theoretic properties of Banach spaces. 
In this article, motivated by the result of Wang and Nan,  we investigate the characterization problem for points of strong subdifferentiability of the norm of $L^p(\mu, X)$.

Throughout this article, all Banach spaces are assumed to be over the real field $\mathbb{R}$. We identify every Banach space $X$, under the canonical embedding, as a subspace of $X^{**}$. Let $B_X$ and $S_X$ denote the closed unit ball and  the unit sphere of $X$, respectively. For $x\in X$ and $r>0$, we denote by $B(x, r)$ the closed ball of radius $r$ centred at $x$. For $x\in X$, we denote by  $J_{X^*}(x)=\{ x^* \in S_{X^*} : x^*(x)=\|x\| \}$ the set of all supporting functionals of $B_X$ at $x$.

We recall  that the norm of a Banach space $X$ is said to be  G\^ateaux differentiable at a vector $x \in X$ if  
$
\lim _{t \rightarrow 0} \frac{\|x+t y\|-\|x\|}{t} 
$
exists for every $y \in X$. If this limit is uniform for every $y \in B_X$, then the norm of $X$ is said to be  Fr\'echet differentiable at   $x$. In \cite{FP}, Franchetti and Paya considered the following generalization of Fr\'echet differentiability. 
\begin{defn}[\cite{FP}]
 The norm of a Banach space $X$ is said to be  strongly subdifferentiable  at $x\in X$ 
if the one sided limit,
$
\lim_{t\to 0^+}\frac{\|x+ty\|-\|x\|}{t}
$
exists uniformly for $y\in B_X$. 
\end{defn}
We say that the norm of $X$ is strongly subdifferentiable if it is strongly subdifferentiable at every point of $S_X$. In \cite{FP},  Franchetti and Paya also observed that, for a Banach space $X$ and $x, y\in X$,  $$\lim_{t\to  0^+}\frac{\|x+ty\|-\|x\|}{t}= \sup \{x^*(y): x^*\in J_{X^*}(x)\}.$$
It is well known that the norm of every finite dimensional Banach space is strongly subdifferentiable at every point. Moreover, the norm  of a Banach space $X$ is Fr\'echet differentiable at $x \in X$ if and only if it is both G\^ateaux differentiable and strongly subdifferentiable at $x$ (see \cite{FP}). Furthermore, the norm of a Banach space $X$ is G\^ateaux differentiable at $x\in X$ if and only if $x$ is a smooth point of  $B_X$, that is,  $x$ admits a unique supporting functional of $B_X$ at $x$. On the other hand, the norm of a Banach space is Fr\' echet differentiable at a nonzero point $x$ if and only if this unique supporting functional is strongly exposed by $x$. The situation, however, is different for strong subdifferentiability, since the supporting functionals at a point of strong subdifferentiability need not be unique. In fact, the norm of  $X$ is strongly subdiffrentiable at  $x\in X$ if and only if   $J_{X^*}(x)$ is strongly exposed by $x$ (see \cite{FP}). We recall this characterization in the following lemma. 
\begin{lem}[\cite{FP}*{Theorem~1.2}]
\label{SSDD}
 Let $X$ be a Banach space and $x\in S_X$. Then the following are equivalent$:$
 \begin{enumerate}[label=(\roman*)]
  \item The norm of $X$ is strongly subdifferentiable at $x$.
  \item For any sequence $(x_n^*)$ in $B_{X^*}$ with $x_n^*(x)\to 1$, the distance $d(x_n^*, J_{X^*}(x))\to 0$. 

 \end{enumerate}
\end{lem}
The characterization in Lemma~\ref{SSDD} will be our main tool for
studying strong subdifferentiability of the norm of Lebesque-Bochner spaces.
We first fix the relevant notation. 

Let $(\Omega, \Sigma, \mu)$ be a finite  measure space and let $X$ be a Banach space. For $1 \leq p<\infty$, let $L^p(\mu, X)$ denote  the Banach space of equivalence classes of   $X$-valued Bochner $p$-integrable functions $f$ defined on $\Omega$, equipped with the norm 
$$\|f\|_p=\left(\int_{\Omega}\|f(t)\|^p d \mu\right)^{1 / p} \textnormal{ for } f \in L^p(\mu, X).$$
We denote by  $L^{\infty}(\mu, X)$  the Banach space of equivalence classes of $X$-valued essentially bounded, strongly measurable  functions, equipped with the essential supremum norm $\|f\|_{\infty}= \ess_{t\in \Omega} \|f(t)\|$.
 If $X^*$ has the Radon-Nikod\'ym property (RNP) with respect to $\mu$, then, by \cite{DU}*{Chapter~IV, Section~1, Theorem~1}, $L^p(\mu, X)^*$ can be identified isometrically with $L^q\left(\mu, X^*\right)$, where $q$ is the conjugate exponent of $p$, with $q=\infty$ when $p=1$. Under this identification, $g\in L^q(\mu, X^*)$ acts on $f\in L^p(\mu, X)$ through the duality pairing
$$
g(f)=\int_{\Omega} g(t)(f(t)) d \mu(t).
$$
It is well known that the separable dual spaces have the Radon-Nykod\'ym property (see \cite{DU}*{Chapter~III, Section~3, Theorem~1}). For a measurable subset $A$ of $\Omega$, we denote by $\mu|_{A}$ the restriction measure of $\mu$
to $A$. We refer to \cite{DU} for further details concerning Bochner spaces and the Radon-Nikod\'ym property.

In  Section~2, we first obtain a complete characterization of the strongly subdifferentiable points of the norm of \(L^{1}(\mu, X)\) when \(X^*\) has the RNP. For $1<p<\infty$, we characterize the strongly subdifferentiable points of the norm of \(L^{p}(\mu, X)\) whenever \(X^*\) is separable and $\mu$ is the Lebesgue measure on $[0, 1]$.



%

\section{\texorpdfstring{Strong subdifferentiability in \(L^{p}(\mu, X)\)}{Strong subdifferentiability in Lp(mu, X)}}


In \cite{CF}, C. Franchetti proved that if   $\mu$ is a positive nonatomic measure, then the norm of  $L^{1}(\mu)$ is nowhere strongly subdifferentiable. Now to consider the case when $\mu$ has atoms, we recall from  \cite{CF} that the strongly subdifferentiable points of the norm of $\ell^1$ are precisely the finitely supported sequences. Motivated by these results, we begin with a characterization of the points at which the norm of $L^1(\mu)$
is strongly subdifferentiable, where $\mu$ is a finite positive measure. As we have not found this characterization stated explicitly in the literature, we include a direct proof.

For $1\leq p<\infty$ and Banach spaces $X$ and $Y$, $X\oplus_{p}Y$ denotes the
$\ell^{p}$-direct sum of $X$ and $Y$, equipped with the norm
$\|(x, y)\|_{p}=\big(\|x\|^{p}+\|y\|^{p}\big)^{1/p}$.

 \begin{lem}\label{L21}
     Let $f\in L^1(\mu)$ be such that  $\|f\|_{1}=1$. Then the norm of $L^{1}(\mu)$ is strongly subdifferentiable at $f$  if and only if  $f=\sum_{i=1}^{n}\alpha_i\, \chi_{A_{i}}$, where $A_{i}$'s are the atoms of $\mu$ and $\alpha_i$'s are nonzero scalars.
 \end{lem}
\begin{proof}
Suppose that $f\in S_{L^1(\mu)}$ is a point of strong subdifferentiability of the norm of $L^1(\mu)$. Let $\lambda$ be the measure given by $\lambda(E)=\int_E|f|\,d\mu$ for $E\in \Sigma$. Denote $S=\{t\in \Omega: f(t)\ne 0\}$. Clearly, for every measurable set $E\subset S$, $\lambda(E)=0$ if and only if $\mu(E)=0$. There exists finitely many atoms $A_1,\ldots, A_n$ and a null set $N$ such that $S=(\bigcup_{i=1}^n A_i)\bigcup N$. Suppose not. Then there exists a sequence $(A_n)$ of subsets of $S$ such that $\lambda(A_n)>0$ and $\lambda(A_n)\to 0$. For $n\in \mathbb{N}$, let $$g_n(t)=
\begin{cases} 
-\sgn(f(t)) &\textnormal{ for } t \in A_n , \\
\sgn(f(t)) &\textnormal{ for } t \in S\setminus A_n , \\
0 &\textnormal{ otherwise. } 
\end{cases}
$$
Then $g_n(f)\to 1$. However, every $h\in J_{L^{\infty}(\mu)}(f)$ satisfies  $h(t)=\sgn(f(t))$ almost everywhere for $t\in S$ and hence $\|g_n-h\|_\infty=2$. Therefore, $d(g_n, J_{L^{\infty}(\mu)}(f))=2$ for every $n$, which contradicts the strong subdifferentiability of the norm at $f$. Since an integrable function is constant on an atom, we have $f=\sum_{i=1}^{k}\alpha_{i}\chi_{A_{i}} $.

Conversely, suppose  $f=\sum_{i=1}^{n}\alpha_{i}\,\chi_{A_{i}} $, where $A_i$'s are atoms of $\mu$. Let $A=\bigcup_{i=1}^{k}A_{i}$.
Then we can write 
\begin{align}
\label{l1split1}
    L^{1}(\mu)=L^{1}(\mu|_A){\displaystyle{\sideset{}{_1}\bigoplus} L^{1}(\mu|_{A^c})}=\ell_n^1 {\displaystyle{\sideset{}{_1}\bigoplus} L^{1}(\mu|_{A^c})}.
\end{align}
We can consider $f$ as an element of $\ell_n^1$. Then  the converse follows immediately from $(\ref{l1split1})$ and \cite{CRJN}*{Theorem~3.12}. 
\end{proof}

 In the following theorem, we show that the vector-valued version of Lemma~\ref{L21} holds for the space $L^1(\mu, X)$, under the assumption that $X^*$ has the RNP. In particular, we show that if $\mu$ is nonatomic, the norm of $ L^{1}(\mu, X)$ is nowhere strongly subdifferentiable, as in the scalar-valued case. 
 
 For $1\leq p<\infty$ and $f\in L^{p}(\mu, X)$, we define $|f|:\Omega\to \mathbb{R}$ by $|f|(t)=\|f(t)\|$. Clearly, $|f|\in L^{p}(\mu)$.


\begin{thm}\label{R10}
    Let $f\in L^{1}(\mu, X)$ be such that  $\|f\|_{1}=1$ and  let $X^*$ has the RNP. Then the norm of $L^{1}(\mu, X)$ is strongly subdifferentiable at $f$  if and only if  $f=\sum_{i=1}^{n}x_{i}\chi_{A_{i}}$, where $A_{i}$'s are the atoms of $\mu$ and each $x_{i}\in B_{X}$ is a point of strong subdifferentiability of the norm of $X$.
\end{thm}
\begin{proof}
We first establish the result for functions $f\in L^1(\mu, X)$ satisfying $f(t)\ne 0$ almost everywhere.

Suppose that $f$ is a point of strong subdifferentiability of the norm satisfying $f(t)\ne 0$ almost everywhere. We first show that $|f|$ is a point of strong subdifferentiability of the norm of $L^1(\mu)$. Let $g_{n}\in B_{L^{\infty}(\mu
)}$ be such that $\int_{\Omega}\|f(t)\|g_{n}(t)d\mu\to 1$. It suffices to prove that $d(g_{n}, J_{L^{\infty}(\mu)}(|f|))\to 0$. Fix  $\phi\in J_{L^{\infty}(\mu, X^*)}(f)$ and define $\phi_{n}=g_{n} \phi $. Clearly $\phi_{n}\in L^{\infty}(\mu, X^*)$ and $\|\phi_{n}\|_{\infty}\leq 1$ for each $n\in \mathbb{N}$. Since $\phi\in J_{L^{\infty}(\mu, X^*)}(f)$, we have $\int_{\Omega}\phi(t)(f(t))d\mu=1=\int_{\Omega}\|f(t)\|d\mu.$
As $\phi(t)(f(t))\leq\|f(t)\|$, it follows that $\phi(t)(f(t))=\|f(t)\|$ almost everywhere and hence $\|\phi(t)\|=1$ almost everywhere. Therefore $\phi_{n}(f)=\int_{\Omega}\phi(t)(f(t))g_{n}(t)d\mu=\int_{\Omega}\|f(t)\|g_{n}(t)d\mu$ and hence, $\phi_{n}(f)\to 1$. 
Since $f$ is a point of strong subdifferentiability of norm, it follows that $d(\phi_{n}, J_{L^{\infty}(\mu, X^*)}(f))\to 0$. 

We now claim that $d(g_{n}, J_{L^{\infty}(\mu)}(|f|))\leq d(\phi_{n}, J_{L^{\infty}(\mu, X^*)}(f))$. Let $\varphi\in J_{L^{\infty}(\mu, X^*)}(f)$. Since $\|\phi(t)\|=1$ almost everywhere, we have 
\begin{align*}
    \|\phi_{n}-\varphi\|_{\infty}&=\|g_{n}\phi-\varphi\|_{\infty}\\
    &\geq \|g_{n}(t)\phi(t)-\varphi(t)\|\\
    &\geq \|g_{n}(t)\phi(t)\|-\|\varphi(t)\|\\
    &=|g_{n}(t)|\|\phi(t)\|-\|\varphi(t)\|\\
    &\geq g_{n}(t)\|\phi(t)\|-|\varphi|(t)\\
    &=g_{n}(t)-|\varphi|(t)
\end{align*}
 almost everywhere. Therefore $\|g_{n}-|\varphi|\|_{\infty}\leq \|\phi_{n}-\varphi\|_{\infty}$. As $|\varphi|\in J_{L^{\infty}(\mu)}(|f|)$, it follows that $$d(g_{n}, J_{L^{\infty}(\mu)}(|f|))\leq \|g_{n}-|\varphi|\|_{\infty}\leq \|\phi_{n}-\varphi\|_{\infty}.$$
Therefore $d(g_{n}, J_{L^{\infty}(\mu)}(|f|))\leq d(\phi_{n}, J_{L^{\infty}(\mu, X^*)}(f))$. Thus $d(g_{n}, J_{L^{\infty}(\mu)}(|f|))\to 0$. Hence $|f|$ is a point of strong subdifferentiability of the norm of $L^{1}(\mu)$. Therefore by Lemma~\ref{L21}, there exist finitely many atoms $A_{1}, A_{2}, \dots, A_{k}$ and scalars $\lambda_{1}, \lambda_{2}, \dots, \lambda_{k}$ such that $|f|=\sum_{i=1}^{k}\lambda_{i}\chi_{A_{i}}$. Since an integrable function is constant on an atom, we have $f=\sum_{i=1}^{k}x_{i}\chi_{A_{i}} $, where $\|x_{i}\|=\lambda_{i}$ for each $i$. Let $A=\bigcup_{i=1}^{k}A_{i}$.
Then we can write 
\begin{align}\label{l1split2}
    L^{1}(\mu, X)=L^{1}(\mu|_A, X){\displaystyle{\sideset{}{_1}\bigoplus} L^{1}(\mu|_{A^c}, X)}=\left(\displaystyle{\sideset{}{_1}\bigoplus_{i=1}^{k} X}\right) {\displaystyle{\sideset{}{_1}\bigoplus} L^{1}(\mu|_{A^c}, X)}.
\end{align}
We can consider $f$ as an element of $\displaystyle{\sideset{}{_1}\bigoplus_{i=1}^{k} X}$. Hence, by \cite{CRJN}*{Theorem~3.12}, each $x_{i}$ is a point of strong subdifferentiability. The converse follows immediately from $(\ref{l1split2})$ and \cite{CRJN}*{Theorem~3.12} and \cite{FP}*{Proposition~2.1}. 
\end{proof}
We now discuss the strong subdifferentiability in $L^{p}(\mu, X)$, when $1<p<\infty$ and $\mu$ is the Lebesgue measure on $[0, 1]$. 
We first prove a sequence of lemmas which will be used in the proof of the main theorem. 

Let $([0,1],\mathcal{A},\mu)$ be the Lebesgue measure space, where $\mathcal{A}$ is the Lebesgue $\sigma$-field and $\mu$ is the Lebesgue measure. Since the measure space $([0,1],\mathcal{A},\mu)$ is complete, we identify functions in the corresponding $L^{p}$-space with representatives that are defined everywhere. Moreover, since the equivalent characterizations of strong subdifferentiability (see Lemma~\ref{SSDD}) involve only sequences, this assumption does not cause any  loss of generality. For a Banach space $X$, the sets $B_{X}$ and $S_{X}$ are always equipped with the Borel $\sigma$-field induced by the norm topology.

\begin{lem}\label{F115}
    Let $g\in S_{L^{p}(\mu)}$ and $(f_n)$ be a sequence in $ B_{L^{q}(\mu)}$  such that $\int_{[0, 1]} |g_{}(t) f_{n}(t)|\,d\mu\to 1$. Then $$\Big\||g|^{p-1}_{}-|f_{n}|\Big\|_{q}\to 0.$$
\end{lem}
\begin{proof}
Suppose $g\in S_{L^{p}(\mu)}$ and $(f_n)$ is a sequence in $ B_{L^{q}(\mu)}$ such that $\int_{[0, 1]} |g_{}(t) f_{n}(t)|\,d\mu\to 1$. We define $\Phi:L^{q}(\mu)\to \mathbb{R}$ by $\Phi(h)=\int_{[0, 1]}h(t)|g(t)|d\mu$. It is easy to check that $\Phi$ is a bounded linear functional on $L^{q}(\mu)$ with $\|\Phi\|=1$ and $\Phi(|g|^{p-1})=1$. By the hypothesis, we have $\Phi(|f_{n}|)\to 1$. Also,
\[
\Big|\Phi\big(\frac{|f_{n}|+|g|^{p-1}}{2}\big)\Big|\leq \Big\|\frac{|f_{n}|+|g|^{p-1}}{2}\Big\|_{q}\leq 1.
\]
Since the left hand side of the above inequality converges to $1$, it follows that $\big\|\frac{|f_{n}|+|g|^{p-1}}{2}\big\|_{q}\to 1$. By the uniform convexity of $L^{q}(\mu)$, we get that $\big\||g|^{p-1}_{}-|f_{n}|\big\|_{q}\to 0.$
\end{proof}

Let $S$ be a Polish space and let $(\Omega, \Sigma, \mu)$ be a measure space. We recall from \cite{SR}*{Chapter~5, Section~5.1} that a set-valued map $F:\Omega\to 2^{S}$ is said to be measurable if for each open subset $U$ of $S$, the set $F^{-1}(U)=\{t\in\Omega: F(t)\cap U\neq \emptyset\}$ is a measurable subset of $\Omega$. 

For an $f\in L^{p}(\mu, X)$,  the following lemma discusses the measurability of the map $F:\Omega\to 2^{B_{X^*}}$ defined by $F(t)=J_{X^*}(f(t))$. 


\begin{lem}\label{L115}
Let \(X\) be a Banach space such that \(X^*\) is separable. Let
\(f\in B_{L^p(\mu,X)}\) be such that $\|f(t)\|=1$ almost everywhere. 
Then the set-valued map $F:[0,1]\to 2^{B_{X^*}}$ defined by $F(t)=J_{X^*}(f(t))$ is measurable.
\end{lem}

\begin{proof}
Without loss of generality, we  assume that $\|f(t)\|=1$ for every \(t\in[0,1]\). Let \(U\subseteq B_{X^*}\) be a non-empty open set. We show that $F^{-1}(U):=\{t\in[0,1]:F(t)\cap U\neq\emptyset\}$
is measurable. Choose an open set \(V\subseteq X^*\) such that $U=V\cap B_{X^*}.$ Let $A_U=\{x\in S_{X}:J_{X^*}(x)\cap U\neq\emptyset\}.$ Then $F^{-1}(U)=\{t\in[0,1]:f(t)\in A_U\}=f^{-1}(A_{U}).$ Since $f$ is measurable, it is enough to prove that \(A_U\) is a measurable subset of \(S_{X}\). Since \(X^*\) is separable, choose a countable  dense subset $D=\{d_1^*,d_2^*,\ldots\}$ of \(X^*\). Consider all sets of the form $C=B_{X^*}\cap B(d_m^*,r)$ for those \(d_{m}^*\in D\) and \(r\in\mathbb Q^+\) such that $B(d_m^*,r)\subseteq V.$ This is a countable family and we enumerate it as $(C_n)_{n=1}^{\infty}.$ 

We now claim that $V\cap B_{X^*}=\bigcup_{n=1}^{\infty}C_n.$ Clearly, $\bigcup_{n=1}^{\infty}C_n\subseteq V\cap B_{X^*}$. Now let \(y^*\in V\cap B_{X^*}\). Since \(V\) is  open, there exists \(\varepsilon>0\) such that $B(y^*,\varepsilon)\subseteq V.$ Choose \(d_m^*\in D\) such that $\|d_m^*-y^*\|<\frac{\varepsilon}{4}.$
Choose \(r\in\mathbb Q^+\) such that $\|d_m^*-y^*\|<r<\frac{\varepsilon}{2}.$ Then $y^*\in B(d_m^*,r)$ and $B(d_m^*,r)\subseteq V.$ 
Therefore $B_{X^*}\cap B(d_m^*,r)=C_{n_{0}}$ for some \(n_{0}\in\mathbb{N}\) and hence \(y^*\in C_{n_{0}}\), which establishes the claim. 

For each \(n\in\mathbb N\), let
\[
A_n=\{x\in S_{X}: \exists\, x^*\in C_n
\text{ such that } x^*(x)=1\}.
\]
It is easy to see that $A_U=\bigcup_{n=1}^{\infty}A_n.$ Now we show that $A_{U}$ is measurable by showing that each \(A_n\) is closed in \(S_{X}\). 
Now let $(x_k)$ be a sequence in $A_{n}$ such that $x_k\to x$. For each \(k\), choose \(x_k^*\in C_n\) such that $x_k^*(x_k)=1.$ Since $X$ is separable and \(C_n\) is weak$^*$-compact, there exists a subsequence
\((x_{k_l}^*)\) of \((x_k^*)\) and \(x^*\in C_n\) such that $x_{k_l}^*\to x^*$ in the weak$^*$-topology. Since $x_{k_{l}}\to x$ in norm,
it follows that $ x^*(x)=\lim x_{k_{l}}^*(x_{k_{l}})=1$. Hence \(x\in A_n\) and \(A_n\) is norm closed in \(S_{X}\). Therefore \(A_U\) is a measurable set in \(S_{X}\), which completes the proof.
\end{proof}

\begin{rmk}\label{R24}
    From the proof of Lemma~\ref{L115}, it follows that 
    the map from $S_{X}$ to $2^{B_{X^*}}$ defined by $x\to J_{X^*}(x)$ is measurable. Therefore by \cite{CV}*{Chapter~III, Theorem~III.9}, we get that for each fixed $x^*\in B_{X^*}$, the function from $S_X$ to $\mathbb{R}$ defined by $x\to d(x^*, J_{X^*}(x))$ is measurable.
\end{rmk} 
We make use of Remark~\ref{R24} to establish the following lemma.
\begin{lem}\label{L118}
    Let $X$ be a Banach space such that $X^*$ is separable and let $h\in L^{q}(\mu, X^*)$ and $f\in L^{p}(\mu, X)$ be such that $h(t)$ and $f(t)$ are non-zero for every $t\in [0, 1]$. Then the function $\tilde{d}:[0, 1]\to \mathbb{R}$ defined by 
    $$
    \tilde{d}(t) = d\left(\frac{h(t)}{\|h(t)\|}, J_{X^*}\left(\frac{f(t)}{\|f(t)\|}\right)\right)\quad \textnormal{ for }t\in [0, 1]
    $$
    is measurable.
\end{lem}

\begin{proof}
    Let $\Pi:[0, 1]\to S_{X^*}\times S_{X}$ be defined by $\Pi(t)=\left(\frac{h(t)}{\|h(t)\|}, \frac{f(t)}{\|f(t)\|}\right)$ and let $G:S_{X^*}\times S_{X}\to \mathbb{R}$ be defined by $G(x^*, x)=d(x^*, J_{X^*}(x))$. Clearly, $\tilde{d}=G\circ \Pi$. Since $h$ and $f$ are both measurable, it is enough to show that $G$ is Borel measurable. 

    Let $D$ be a countable dense subset of $S_{X^*}$. For $r>0$, we let $G_{r}=\{(x^*, x): d(x^*, J_{X^*}(x))\geq r\}$. Then it is easy to see that 
    $$
    G_{r}=\bigcap_{n\in \mathbb{N}}\bigcup_{z^*\in D}\left\{\Big\{x^*\in S_{X^*}:\|x^*-z^*\|<\frac{1}{n}\Big\}\times\Big\{x\in S_{X}: d(z^*, J_{X^*}(x))>r-\frac{1}{n}\Big\}\right\}
    $$
     Then it follows from Remark~\ref{R24} that 
    $G_{r}$ is a Borel set. Thus $G$ is Borel measurable and hence $\tilde{d}$ is measurable.
\end{proof}
\begin{rmk} We observe that if $h$ and $f$ are non-zero almost everywhere, then the function $\tilde{d}$ is defined almost everywhere and measurable.
\end{rmk}

\begin{lem}\label{L116}
Let \(X\) be a Banach space such that \(X^*\) is separable, and let
\(f\in B_{L^p(\mu,X)}\) be such that $\|f(t)\|=1$  almost everywhere. Then, for every \(\varepsilon>0\) and
\(\delta>0\), the set
\[
Z=\left\{
t\in[0,1]: \textnormal{ there exists }\, x^*\in S_{X^*} \textnormal{ such that }
x^*(f(t))>1-\delta
\textnormal{ and }
d\bigl(x^*,J_{X^{*}}(f(t))\bigr)\geq \varepsilon
\right\}
\]
is a measurable set. 
\end{lem}
\begin{proof}
   Let $\{x_{k}^*\}$ be a countable dense subset of $S_{X^*}$ and for $k\in \mathbb{N}$, let $$Z_{k}=\left\{t\in[0,1]:
x_{k}^*(f(t))>1-\delta
\textnormal{ and }
d\bigl(x_{k}^*,J_{X^*}(f(t))\bigr)\geq \varepsilon
\right\}.$$ 
Since, for each fixed $t\in Z$, the set $\{x^*\in S_{X^*}: x^*(f(t))>1-\delta
\textnormal{ and }
d\bigl(x^*,J_{X^{*}}(f(t))\bigr)> \varepsilon\}$ is an open set, it follows that $Z=\bigcup_{k=1}^{\infty}Z_{k}.$ So, it is enough to show that each $Z_{k}$ is a measurable set. Now, for each fixed $k\in\mathbb{N}$, it follows from \cite{CV}*{Chapter~III, Theorem~III.9} and Lemma~\ref{L115} that the function
 \[
 t\to d\bigl(x_{k}^*,J_{X^*}(f(t))\bigr)
 \]
is measurable. Consequently, the set $\left\{t\in[0,1]: d\bigl(x_{k}^*,J_{X^*}(f(t))\bigr)\geq \varepsilon
\right\}$ is a measurable set. Since, the set $\left\{t\in[0,1]: x_{k}^*(f(t))>1-\delta \right\}$ is measurable, it follows that $Z_{k}$ is a measurable set for each $k\in\mathbb{N}$, which completes the proof.
\end{proof}

For a Polish space $S$, we denote by $\mathcal{B}_S$, the Borel $\sigma$-algebra and by $\mathcal{A}\otimes \mathcal{B}_S$, the product $\sigma$-algebra on $[0,1]\times S$. Let $F:[0,1]\to 2^{S}$ be a set-valued 
mapping. The graph of $F$ is the subset of $[0,1]\times S$ defined by $\operatorname{Graph}(F)
=\{(t,s)\in[0,1]\times S:s\in F(t)\}$ (see \cite{SR} for details). We now present our main result.


\begin{thm}\label{T28}
  Let \(X\) be a Banach space such that \(X^*\) is separable, and let
\(f\in S_{L^p(\mu,X)}\) be such that $f(t)\neq 0$  almost everywhere. Then $f$ is a point of strong subdifferentiability of the norm of $L^p(\mu,X)$ if and only if  $\frac{f(t)}{\|f(t)\|}$ is a point of strong subdifferentiability of the norm of $X$ for almost every \(t\in[0,1]\) with $f(t)\neq 0$.
\end{thm}

\begin{proof}
Let $f\in S_{L^p(\mu,X)}$be such that $f(t)\neq 0$  almost everywhere. Without loss of generality, we assume that $f(t)\neq 0$ for every $t\in [0, 1]$. Let
$$
A=\Big\{t\in [0, 1]:\textnormal{ norm of } X \textnormal{ is not strongly subdifferentiable at }\frac{f(t)}{\|f(t)\|}\Big\}
.$$ We first show that $A$ is a measurable set. For each $t\in A$, since $\frac{f(t)}{\|f(t)\|}$ is not a point of strong subdifferentiability, there exists  $n\in \mathbb{N}$ such that for each $m\in \mathbb{N}$, there exists $x_{}^*\in B_{X^*}$ satisfying 
\begin{align}\label{ES}
x^*\left(\frac{f(t)}{\|f(t)\|}\right)>1-\frac{1}{m}
\quad \textnormal{and} \quad
d\left(
x^*,
J_{X^*}\left(\frac{f(t)}{\|f(t)\|}\right)
\right)\geq\frac{1}{n}.
\end{align}
For $n, m\in\mathbb{N}$, let
$$
A_{n, m}=\{t\in [0, 1]: \textnormal{there exists } x^*\in B_{X^*} \textnormal{ satisfying }(\ref{ES})\}.
$$
It is easy to see that $A=\bigcup_{n=1}^{\infty}\bigcap_{m=1}^{\infty}A_{n, m}.$ From Lemma~\ref{L116}, we have that each $A_{n, m}$ is a measurable set. Consequently, $A$ is a measurable set. Suppose that $\mu(A)>0$. We now show that $f$ is not a point of strong subdifferentiability by exhibiting a sequence $(\tilde{h}_{m})$ in $ B_{L^{q}(\mu, X^*)}$ such that $\tilde{h}_{m}(f)\to 1$ and $d\big(\tilde{h}_{m}, J_{L^{q}(\mu, X^*)}(f)\big)$ does not converges to $0$.

Since $\mu(A)>0$, we have that $\mu(\bigcap_{m=1}^{\infty}A_{n_{0}, m})>0$ for some $n_{0}\in \mathbb{N}$. Let $\varepsilon=\frac{1}{n_{0}}$ and $B_{m}= A_{n_{0}, m}$. Clearly, $\mu(B_{m})>0$ for each $m\in \mathbb{N}$. For $m\in\mathbb{N}$ and $t\in B_{m}$, let $$B_{m}(t)=\left\{x^*\in B_{X^*}: x^*\left(\frac{f(t)}{\|f(t)\|}\right)>1-\frac{1}{m} \textnormal{ and }d\left(
x^*,
J_{X^*}\left(\frac{f(t)}{\|f(t)\|}\right)
\right)\geq\varepsilon\right\}.$$
Clearly, for each $t\in B_{m}$, the set $B_{m}(t)$ is non-empty. We can write 
$$B_{m}(t)=\left\{x^*\in B_{X^*}: x^*\left(\frac{f(t)}{\|f(t)\|}\right)>1-\frac{1}{m}\right\}\bigcap \left\{x^*\in B_{X^*}: d\left(
x^*,
J_{X^*}\left(\frac{f(t)}{\|f(t)\|}\right)
\right)\geq\varepsilon\right\}.$$
Since for each $t\in [0, 1]$, the maps $x^*\to x^*\left(\frac{f(t)}{\|f(t)\|}\right)$ and $x^*\to d\left(
x^*,
J_{X^*}\left(\frac{f(t)}{\|f(t)\|}\right)
\right)$ are both measurable, we get that set $B_{m}(t)$ is a measurable subset of $B_{X^*}$.


Define $ F_m: [0,1] \to 2^{B_{X^*}} $ by
\[
F_{m}(t) =
\begin{cases} 
B_{m}(t) &\textnormal{ for } t \in B_{m},  \\
J_{X^*}\left(\frac{f(t)}{\|f(t)\|}\right) &\textnormal{ for } t \notin B_{m}.
\end{cases}
\]
Let 
\begin{align*}
        E_{1}&=\big\{(t, x^*): t \in B_{m}\textnormal{ and }x^*\in B_{m}(t) \big\}\textnormal{ and } E_{2}=\left\{(t, x^*): t \notin B_{m}\textnormal{ and }x^*\in J_{X^*}\left(\frac{f(t)}{\|f(t)\|}\right) \right\} .
\end{align*}
Clearly, $\textnormal{Graph}(F_m) = E_{1}\bigcup E_{2}.$

We now show that $\textnormal{Graph}(F_{m})$ is a measurable set. It follows from the proof of \cite{SD1}*{Theorem~2.3} that the function $(t, x^*)\to x^*\left(\frac{f(t)}{\|f(t)\|}\right)$ from $[0, 1]\times B_{X^*} \to \mathbb{R}$ is measurable. Since $B_{m}^c$ is a measurable set, we get that $E_{2}$ is a measurable set. We next show that $E_{1}$ is a measurable set. Let 
\begin{align*}
    E_{1}'&=\left\{(t, x^*)\in [0, 1]\times B_{X^*}: t \in B_{m} \textnormal{ and }x^*\left(\frac{f(t)}{\|f(t)\|}\right)>1-\frac{1}{m} 
\right\}\ \textnormal{ and }\\
    E_{1}''&=\left\{(t, x^*)\in [0, 1]\times B_{X^*}: t \in B_{m} \textnormal{ and } d\left(
x^*,
J_{X^*}\left(\frac{f(t)}{\|f(t)\|}\right)
\right)\geq\varepsilon\right\}.
\end{align*}
Clearly, 
$
E_{1}=E_{1}'\bigcap E_{1}''. 
$
Using an argument similar to the one used to establish the measurability of $E_{2}$, we can see that the set $E_{1}'$
is a measurable set. Now to show $E_1''$ is a measurable set, it is enough to show that the mapping $\phi:[0, 1]\times B_{X^*}\to \mathbb{R}$ defined by
$$
(t, x^*)\to d\left(
x^*,
J_{X^*}\left(\frac{f(t)}{\|f(t)\|}\right)
\right)$$ is measurable. It is easy to check that the map $\phi$ is continuous in the variable $x^*$ and it follows  from Lemma~\ref{L115} and \cite{CV}*{Chapter~3, Theorem~III.9} that the map $\phi$ is measurable in the variable $t$. Then, by \cite{CV}*{Lemma~III.14}, we get that $\phi$ is measurable. Thus, $E_{1}$ is a measurable set and therefore $\operatorname{Graph}(F_{m})$ is a measurable set. 

Since $B_{X^*}$ is a Polish space with respect to norm topology, by applying Von Neumanns's selection theorem (see \cite{SR}*{Chapter~5, Corollary~5.5.8}) to the measurable set $\operatorname{Graph}(F_m)\subseteq [0, 1]\times B_{X^*}$, we obtain a measurable selection
$h_m : [0,1] \to B_{X^*}$ such that 
\[
\begin{aligned}
\text{(a)}\;& \|h_m(t)\| \leq 1 
&& \text{for all } t \in [0,1],\\
\text{(b)}\;& h_{m}(t)\left(\frac{f(t)}{\|f(t)\|}\right) = 1 
&& \text{for every } t \notin B_{m},\\
\text{(c)}\;& h_{m}(t)\left(\frac{f(t)}{\|f(t)\|}\right) > 1 - \frac{1}{m} \textnormal{ and }d\left(
h_{m}(t),
J_{X^*}\left(\frac{f(t)}{\|f(t)\|}\right)
\right)\geq\varepsilon
&& \text{for all } t \in B_m,\\
\end{aligned}
\]
Since $X^*$ is separable, by Pettis measurability theorem (see \cite{DU}*{Chapter~II, Section~1, Theorem~2}), it follows that the function $h_{m}:[0,1]\to B_{X^*}$ is strongly measurable. Hence $h_{m}\in B_{L^{q}(\mu, X^*)}$ for all $m\in \mathbb{N}$. 

Now for $m\in \mathbb{N}$,  let $\tilde{h}_{m}=|f|^{p-1} h_{m}$.

\noindent \textbf{Claim: } $\tilde{h}_{m}\in B_{L^{q}(\mu, X^*)}$ for all $m\in \mathbb{N}$, 
 $
 \tilde{h}_{m}(f)\to1\ \textnormal{and there exists } C>0 \textnormal{ such that}\ d\big(\tilde{h}_{m}, J_{L^{q}(\mu, X^*)}(f)\big)\geq C
 $
 for all $m\in \mathbb{N}$. 
 
 This claim contradicts the fact that $f$ is a point of strong subdifferentiability of the norm of $L^{p}(\mu, X)$. 
 
 Now to establish the claim, we first observe that $\|\tilde{h}_{m}\|_{q}\leq 1$. Then it follows that
 $$
1\geq \tilde{h}_{m}(f) = \int_{[0, 1]}\|f(t)\|^{p-1}h_{m}(t)(f(t))d\mu\geq\int_{B_{m}^c}\|f(t)\|^{p}d\mu+\big(1-\frac{1}{m}\big)\int_{B_{m}}\|f(t)\|^{p}\geq (1-\frac{1}{m}) \|f\|^{p}_{p}.
 $$
 Since $\|f\|_{p}=1$, we have that $\tilde{h}_{m}(f) \to1.$
 
 Let $g\in J_{L^{q}(\mu, X^*)}(f)$ be an arbitrary element. We have 
 \begin{align}\label{E100}
     1=g(f)=\int_{[0, 1]}g(t)(f(t))d\mu\leq \int_{[0, 1]}\|g(t)\|\|f(t)\|d\mu\leq \|g\|_{q}\|f\|_{p}=1.
 \end{align}
 
 Since equality holds in Hölder's inequality, we get that the functions 
$t \to \|g_{}(t)\|^{q}$ and 
$t \to \|f(t)\|^{p}$ are linearly dependent. So there exists a $\lambda\in\mathbb{R}$ such that $\|g_{}(t)\|^{q}=\lambda\|f(t)\|^{p}$ almost everywhere. Since $\|g_{}\|_{q}=\|f\|_{p}=1$, we get that $\lambda=1$. Thus $\|g_{}(t)\|^{q}=\|f(t)\|^{p}$ almost everywhere. From~(\ref{E100}), it follows that
$g(t)(f(t))= \|g(t)\|\|f(t)\|$ almost everywhere.

Let $\tilde{g}=\frac{g}{|f|^{p-1}}$. It is easy to see that $\tilde{g}\in J_{L^{q}(\mu, X^*)}\big(\frac{f}{|f|}\big)$ and $\tilde{g}(t)\in J_{X^*}(\frac{f(t)}{\|f(t)\|})$ almost everywhere. Then, we have 
\begin{align*}
    \|\tilde{h}_{m}-g\|_{q}&=\left\||f|^{p-1} h_{m}-|f|^{p-1} \tilde{g}\right\|_{q}\\
    &=\left(\int_{[0, 1]}\|f(t)\|^{q(p-1)}\|h_{m}(t)-\tilde{g}(t)\|^{q}d\mu\right)^{\frac{1}{q}}\\
    &\geq \left(\int_{\bigcap_{m=1}^{\infty} A_{n_{0}, m}}\|f(t)\|^{q(p-1)}\|h_{m}(t)-\tilde{g}(t)\|^{q}d\mu\right)^{\frac{1}{q}}\\
    &\geq \left(\int_{\bigcap_{m=1}^{\infty} A_{n_{0}, m}}\|f(t)\|^{p}\,d\Big(
h_{m}(t),
J_{X^*}\Big(\frac{f(t)}{\|f(t)\|}\Big)
\Big)^{q}d\mu\right)^{\frac{1}{q}}\\
&\geq C, 
\end{align*}
where $C=\frac{1}{n_{0}}\left(\int_{\bigcap_{m=1}^{\infty} A_{n_{0}, m}}\|f(t)\|^{p}d\mu\right)^{\frac{1}{q}}$. Hence it follows that 
 $$
 d\big(\tilde{h}_{m}, J_{L^{q}(\mu, X^*)}(f)\big)\geq C\quad \textnormal{ for all } m\in \mathbb{N}.
 $$
 This proves the claim and hence completes the proof of the necessary part.

 Conversely, suppose $\frac{f(t)}{\|f(t)\|}$ is a point of strong subdifferentiability for almost all $t\in[0, 1]$. Let $(h_{n})$  be a sequence in $B_{L^{q}(\mu, X^*)}$ such that $h_{n}(f)\to 1$. To show that $d(h_{n}, J_{L^{q}(\mu, X^*)}(f))\to 0$, we exhibit a sequence $(\tilde{g}_{k})$ in $J_{L^{q}(\mu, X^*)}(f)$ such that $\|h_{n}-\tilde{g}_{k}\|_{q}\to 0$.

 Since
 $$h_{n}(f)=\int_{[0, 1]}h_{n}(t)(f(t))d\mu\leq \int_{[0, 1]}\|h_{n}(t)\|\|f(t)\|d\mu\leq\Big(\int_{[0, 1]}\|f(t)\|^{p}d\mu\Big)^{\frac{1}{p}}=1,$$
 we have that 
 \begin{align}\label{E23}
    \lim_{n\to\infty}\int_{[0, 1]}h_{n}(t)(f(t))d\mu= \lim_{n\to\infty}\int_{[0, 1]}\|h_{n}(t)\|\|f(t)\|d\mu=1.
 \end{align}
Now, by applying Lemma~\ref{F115} to (\ref{E23}), we get that 
\begin{align}\label{E25}
    \||h_{n}|-|f|^{p-1}\|_{q}\to 0.
\end{align}
Without loss of generality, we may assume that $h_{n}(f)\neq 0$ for each $n\in\mathbb{N}$ and that $h_{n}(t)\neq 0$ almost everywhere for each $n\in \mathbb{N}$. For, suppose $Z_{n}=\{t\in[0, 1]:h_{n}(t)=0\}$ and $0<\mu(Z_{n})<1$. Fix $h\in J_{L^{q}(\mu, X^*)}\big(\frac{f}{|f|}\big)$. Then $h(t)(f(t))=\|f(t)\|$ and $\|h_{}(t)\|=1$ almost everywhere. We define another sequence $(\hat{h}_{n})$ by 
$$\hat{h}_{n}=h_{n}\chi_{Z_{n}^{c}} +|f|^{p-1} h\chi_{Z_{n}}$$ 
    Clearly, $\hat{h}_{n}(t)\neq 0$ almost everywhere for each $n\in \mathbb{N}$. It follows from~(\ref{E25}) that $\int_{Z_{n}}\|f(t)\|^{p}\,d\mu\to 0$. This implies that $\hat{h}_{n}(f)\to 1$ and $\|\hat{h}_{n}-h_{n}\|_{q}\to 0$. Put $\alpha_{n}=\textnormal{max}\{1, \|\hat{h}_{n}\|_{q}\}$. Since $\|h_{n}\|_{q}\to 1$, we have $\|\hat{h}_{n}\|_{q}\to 1$ and hence $\alpha_{n}\to 1$. Therefore $\frac{\hat{h}_{n}}{\alpha_{n}}(f)\to 1$, $\|\frac{\hat{h}_{n}}{\alpha_{n}}-h_{n}\|_{q}\to 0$ and $\big\|\frac{\hat{h}_{n}}{\alpha_{n}}\big\|_{q}\leq 1$. Hence $(h_{n})$ may be replaced by $\bigl(\frac{\hat{h}_{n}}{\alpha_{n}}\bigr)$. This justifies the assumption that $h_{n}(t)\neq 0$ almost everywhere for each $n\in \mathbb{N}$. We further assume, without loss of generality, that $h_{n}(t)\neq 0$ for every $t\in [0, 1]$ and for each $n\in \mathbb{N}$.

Now we can write~(\ref{E23}) as follows
\begin{align}\label{E24}
    \int_{[0, 1]}\|h_{n}(t)\|\|(f(t))\|\left(1-\frac{h_{n}(t)}{\|h_{n}(t)\|}\Big(\frac{f(t)}{\|f(t)\|}\Big)\right)d\mu\to0.
\end{align}

Also, by Hölder's inequality, we get that

\begin{align}\label{E26}
  \int_{[0, 1]}(\|h_{n}(t)\|-\|f(t)\|^{p-1})\|f(t)\|\left(1-\frac{h_{n}(t)}{\|h_{n}(t)\|}\Big(\frac{f(t)}{\|f(t)\|}\Big)\right)d\mu\leq M\big\||h_{n}|-|f|^{p-1}\big\|_{q},
    \end{align}
for some $M>0$. Therefore, by adding (\ref{E24}) and (\ref{E26}), we have 
\begin{align}\label{E27}
    \int_{[0, 1]}\|f(t)\|^{p}\left(1-\frac{h_{n}(t)}{\|h_{n}(t)\|}\Big(\frac{f(t)}{\|f(t)\|}\Big)\right)d\mu\to 0.
\end{align}
Since the integrand in (\ref{E27}) is a positive function, we can choose a subsequence $(h_{n_{k}})$ of $(h_{n})$ such that 
$$
\frac{h_{n_{k}}(t)}{\|h_{n_{k}}(t)\|}\left(\frac{f(t)}{\|f(t)\|}\right)\to1\ \textnormal{almost everywhere.}
$$
Since $\frac{f(t)}{\|f(t)\|}$ is a point of strong subdifferentiability of the norm of $X$ for almost every $t\in[0, 1]$, we get that 
\begin{align*}
    d\left(
\frac{h_{n_{k}}(t)}{\|h_{n_{k}}(t)\|},
J_{X^*}\Big(\frac{f(t)}{\|f(t)\|}\Big)
\right)\to 0 \quad \textnormal{almost everywhere}.
\end{align*}
We may assume that 
\begin{align}\label{E28}
    d\left(
\frac{h_{n_{k}}(t)}{\|h_{n_{k}}(t)\|},
J_{X^*}\Big(\frac{f(t)}{\|f(t)\|}\Big)
\right)\to 0 \quad \textnormal{for each } t\in [0, 1].
\end{align}
We now construct a sequence $(\tilde{g}_{k})$ in $J_{L^{q}(\mu, X^*)}(f)$ and a subsequence of $(h_{n_{k}})$, which we again denote by $(h_{n_{k}})$, such that $\|h_{n_{k}}-\tilde{g}_{k}\|_{q}\to 0$.

Let $k\in\mathbb{N}$ be fixed. For each $t\in[0, 1]$, we define 
$$
A_{t}=\left\{x^*\in J_{X^*}\big(\frac{f(t)}{\|f(t)\|}\big):\Big\|\frac{h_{n_{k}}(t)}{\|h_{n_{k}}(t)\|}-x^*\Big\|<d\Big(
\frac{h_{n_{k}}(t)}{\|h_{n_{k}}(t)\|},
J_{X^*}\big(\frac{f(t)}{\|f(t)\|}\big)
\Big)+\frac{1}{k}\right\}.
$$
It is easy to see that $A_{t}$ is a Borel set in $S_{X^*}$ for each $t\in [0, 1]$. Now define $F_{k}:[0, 1]\to 2^{S_{X^*}}$ by $F_{k}(t) = A_{t}$. Then 

\noindent$
{\textnormal{Graph}(F_{k})=\left\{(t, x^*):t\in [0, 1],\, x^*\in J_{X^*}\big(\frac{f(t)}{\|f(t)\|}\big)\ \textnormal{and}\ \Big\|\frac{h_{n_{k}}(t)}{\|h_{n_{k}}(t)\|}-x^*\Big\|<d\Big(
\frac{h_{n_{k}}(t)}{\|h_{n_{k}}(t)\|},
J_{X^*}\big(\frac{f(t)}{\|f(t)\|}\big)
\Big)+\frac{1}{k}\right\}}.
$

Write
$$
\varphi_{k}(t)=d\Big(
\frac{h_{n_{k}}(t)}{\|h_{n_{k}}(t)\|},
J_{X^*}\big(\frac{f(t)}{\|f(t)\|}\big)
\Big)\textnormal{ and }
\Phi_{k}(t, x^*)=\Big\|\frac{h_{n_{k}}(t)}{\|h_{n_{k}}(t)\|}-x^*\Big\|.
$$
By Lemma~\ref{L118}, $\varphi_{k}$ is measurable. Moreover $\Phi_{k}$ is measurable in the variable $t$ and continuous in the variable $x^*$, hence jointly
measurable on $[0, 1]\times S_{X^*}$ (see \cite{CV}*{Chapter~3, Lemma~III.14}). Since $$\textnormal{Graph}(F_{k})=\Big\{(t, x^*)\in [0, 1]\times S_{X^*}:x^*\big(\frac{f(t)}{\|f(t)\|}\big)=1\Big\}\bigcap \Big\{(t, x^*)\in [0, 1]\times S_{X^*}: \Phi_{k}(t, x^*)-\varphi_{k}(t)<\frac{1}{k}\Big\},$$ it follows that $\textnormal{Graph}(F_{k})$ is a measurable set in $[0, 1]\times S_{X^*}$. Now by Von Neumann's selection theorem (see \cite{SR}*{Chapter~5, Corollary~5.5.8}), we obtain a measurable selection $g_{k}:[0, 1]\to S_{X^*}$ such that $g_{k}(t)\in A_{t}$ for each $t\in [0, 1]$. Then, by Pettis measurability theorem (see \cite{DU}*{Chapter~III, Section~1, Theorem~2}), it follows that the function $g_{k}:[0,1]\to S_{X^*}$ is strongly measurable and hence $g_{k}\in S_{L^{q}(\mu, X^*)}$. Therefore
\begin{align}
    \Big\|\frac{h_{n_{k}}}{|h_{n_{k}}|}-g_{k}\Big\|_{q}^{q}= \int_{[0, 1]}\Big\|\frac{h_{n_{k}}(t)}{\|h_{n_{k}}(t)\|}-g_{k}(t)\Big\|^{q}d\mu\leq\int_{[0, 1]}\left(d\Big(
\frac{h_{n_{k}}(t)}{\|h_{n_{k}}(t)\|},
J_{X^*}\big(\frac{f(t)}{\|f(t)\|}\big)
\Big)+\frac{1}{k}\right)^{q}d\mu.
\end{align}
Then, by using (\ref{E28}) and Dominated convergence theorem, we get that 
$$
\int_{[0, 1]}\left(d\Big(
\frac{h_{n_{k}}(t)}{\|h_{n_{k}}(t)\|},
J_{X^*}\big(\frac{f(t)}{\|f(t)\|}\big)
\Big)+\frac{1}{k}\right)^{q}d\mu\to 0.
$$
Consequently, 
\begin{align}\label{E210}
     \Big\|\frac{h_{n_{k}}}{|h_{n_{k}}|}-g_{k}\Big\|_{q}\to 0.
\end{align}
Thus there exists a subsequence of $(h_{n_{k}})$ and a subsequence of $(g_{{k}})$, which we continue to denote by $(h_{n_{k}})$ and $(g_{{k}})$, such that $\frac{h_{n_{k}}(t)}{\|h_{n_{k}}(t)\|}-g_{k}(t)\to 0$ almost everywhere. Then, by Dominated convergence theorem, it follows that 
\begin{align}\label{E211}
      \Big\||f|^{p-1}\frac{h_{n_{k}}}{|h_{n_{k}}|}-g_{k}\Big\|_{q}\to 0.
\end{align}
Let $\tilde{g}_{k}=|f|^{p-1} g_{k}$. Now it is enough to show that $\tilde{g}_{k}\in J_{L^{q}(\mu, X^*)}(f)$ and $\|h_{n_{k}}-\tilde{g}_{k}\|_{q}\to 0$. Since $g_{k}(t)\in J_{X^*}\left(\frac{f(t)}{\|f(t)\|}\right)$ for each $t\in [0, 1]$, we get $\tilde{g}_{k}(t)(f(t))=\|f(t)\|^{p}$ and hence $\tilde{g}_{k}\in J_{L^{q}(\mu, X^*)}(f)$. Then, we have 
\begin{align*}
    \big\|h_{n_{k}}-\tilde{g}_{k}\big\|_{q}&\leq\Big\|h_{n_{k}}-\frac{h_{n_{k}}}{|h_{n_{k}}|} |f|^{p-1}\Big\|_{q}+\Big\|\frac{h_{n_{k}}}{|h_{n_{k}}|} |f|^{p-1}-|f|^{p-1} g_{k}\Big\|_{q}\\
     &\leq \big\||h_{n_{k}}|-|f|^{p-1}\big\|_{q}+\Big\||f|^{p-1} \left(\frac{h_{n_{k}}}{|h_{n_{k}}|}-g_{k} \right)\Big\|_{q}.
\end{align*}
From (\ref{E25}) and (\ref{E211}), we get that $\|h_{n_{k}}-\tilde{g}_{k}\|_{q}\to 0$. Hence $d\big(h_{n_{k}}, J_{L^{q}(\mu, X^*)}(f)\big)\to 0$. Thus every subsequence of $(d\big(h_{n_{}}, J_{L^{q}(\mu, X^*)}(f)\big))$ has a further subsequence converging to $0$. Hence $d\big(h_{n_{}}, J_{L^{q}(\mu, X^*)}(f)\big)\to0$. This completes the proof.
\end{proof}
Let $f\in S_{L^{p}(\mu, X)}$ and let $Z=\{t\in [0,1]: f(t)=0\}$. 
Then we have the decomposition 
$L^{p}(\mu,X)=L^{p}(\mu|_{Z},X)\displaystyle{\sideset{}{_p}\bigoplus}L^{p}(\mu|_{Z^c},X)$.
Since $f$ vanishes on $Z$, we may assume, without loss of generality, that 
$f\in S_{L^{p}(\mu|_{Z^c},X)}$. By \cite{FP}*{Proposition~2.1}, the norm of 
$L^{p}(\mu,X)$ is strongly subdifferentiable at $f$ if and only if the norm of 
$L^{p}(\mu|_{Z^c},X)$ is strongly subdifferentiable at $f$. As a consequence, we have the following general version of Theorem~\ref{T28}.
\begin{cor}\label{C29}
Let \(X\) be a Banach space such that \(X^*\) is separable, and let
\(f\in S_{L^p(\mu,X)}\). Then $f$ is a point of strong subdifferentiability of the norm of $L^p(\mu,X)$ if and only if  $f(t)$ is a point of strong subdifferentiability of the norm of $X$ for almost every \(t\in[0,1]\) with $f(t)\neq 0$. 
\end{cor}

A Banach space $X$ is said to be an Asplund space if $X^*$ has the RNP (see \cite{RRP} for details). In \cite{FP}*{Theorem~5.1}, the authors have shown that if the norm of $X$ is strongly subdifferentiable, then $X$ is an Asplund space. We use this result to establish the following corollary.
\begin{cor}
     Let $X$ be a separable Banach space. Then, the norm of $L^{p}(\mu, X)$ is strongly subdifferentiable if and only if the norm of $X$ is strongly subdifferentiable.
\end{cor}
\begin{proof}
Let $X$ be a separable space such that the norm of $X$ is strongly subdifferentiable. By \cite{FP}*{Theorem~5.1}, we get that $X$ is an Asplund space. It follows from \cite{RRP}*{Section~2, Corollary~2.15} that $X^*$ is separable. Now it follows from Corollary~\ref{C29} that the norm $L^{p}(\mu, X)$ is strongly subdifferentiable. The converse now follows from the fact that $X$ embeds isometrically in $L^{p}(\mu, X)$ and strong subdifferentiability of the norm is inherited by subspaces.
\end{proof}

\noindent\textbf{Acknowledgments}
 Most of this work was carried out during the second author’s visit to Shiv Nadar University and the first author’s subsequent visit to IIT Palakkad. The first author thanks IIT Palakkad for its hospitality, while the second author thanks Prof. T. S. S. R. K. Rao and Shiv Nadar University, Delhi-NCR for their hospitality. The authors also thank Prof. T. S. S. R. K. Rao for several fruitful discussions. Both authors gratefully acknowledge financial support from KSCSTE, Govt. of Kerala for these visits.


 \noindent\textbf{Competing interests} The authors have no competing interests to declare that are relevant to the content of this article.



\bibliographystyle{amsplain}
\bibliography{bibliography_new}
\end{document}